\documentclass[12pt]{article}                                               

\input{amssym.def}                                                           
\input{amssym.tex}

\begin{document}                                                             
\title{$K$-theory of  cluster $C^*$-algebras}

\author{Igor  ~Nikolaev}


\date{}
 \maketitle


\newtheorem{thm}{Theorem}
\newtheorem{lem}{Lemma}
\newtheorem{dfn}{Definition}
\newtheorem{rmk}{Remark}
\newtheorem{cor}{Corollary}
\newtheorem{cnj}{Conjecture}
\newtheorem{exm}{Example}


\begin{abstract}
It is proved that  the $K_0$-group of a cluster $C^*$-algebra
is isomorphic to the corresponding cluster algebra.  As a 
corollary,   one gets a shorter   proof of the positivity conjecture
for  cluster algebras. 
As an example,  we consider a cluster $C^*$-algebra ${\Bbb A}(1,1)$ 
 coming  from triangulation of an annulus with one marked point on each boundary
 component.

\vspace{3mm}


\vspace{3mm}

\end{abstract}

\section{Introduction}
Cluster algebras are a class of commutative rings introduced by 
[Fomin \& Zelevinsky 2002]  \cite{FoZe1}  having deep roots  in    
hyperbolic  geometry and Teichm\"uller theory  [Williams 2014] \cite{Wil1}.  
Namely,  the {\it cluster algebra}  ${\cal A}(\mathbf{x}, B)$ of rank $n$ 
is a subring of the field  of  rational functions in $n$ variables
depending  on a {\it cluster} of variables  $\mathbf{x}=(x_1,\dots, x_n)$
and a skew-symmetric matrix  $B=(b_{ij})\in M_n(\mathbf{Z})$; 
the pair  $(\mathbf{x}, B)$ is called a {\it seed}.
A new cluster $\mathbf{x}'=(x_1,\dots,x_k',\dots,  x_n)$ and a new
skew-symmetric matrix $B'=(b_{ij}')$ is obtained from 
$(\mathbf{x}, B)$ by the  {\it exchange relations}:
\begin{eqnarray}\label{eq1}
x_kx_k'  &=& \prod_{i=1}^n  x_i^{\max(b_{ik}, 0)} + \prod_{i=1}^n  x_i^{\max(-b_{ik}, 0)},\cr 
b_{ij}' &=& \cases{
-b_{ij}  & \hbox{if}  $i=k$ \hbox{or} $j=k$\cr
b_{ij}+{|b_{ik}|b_{kj}+b_{ik}|b_{kj}|\over 2}  & \hbox{otherwise.}
}
\end{eqnarray}
The seed $(\mathbf{x}', B')$ is said to be a {\it mutation} of $(\mathbf{x}, B)$ in direction $k$,
where $1\le k\le n$;   the algebra  ${\cal A}(\mathbf{x}, B)$ is  generated by cluster  variables $\{x_i\}_{i=1}^{\infty}$
obtained from the initial seed $(\mathbf{x}, B)$ by the iteration of mutations  in all possible
directions $k$.

The {\it Laurent phenomenon} 
proved by  [Fomin \& Zelevinsky 2002]  \cite{FoZe1}  says  that  ${\cal A}(\mathbf{x}, B)\subset \mathbf{Z}[\mathbf{x}^{\pm 1}]$,
where  $\mathbf{Z}[\mathbf{x}^{\pm 1}]$ is the ring of  the Laurent polynomials in  variables $\mathbf{x}=(x_1,\dots,x_n)$
depending on an initial seed $(\mathbf{x}, B)$;  
in other words, each  generator $x_i$  of  algebra ${\cal A}(\mathbf{x}, B)$  can be 
written as a  Laurent polynomial in $n$ variables with    integer coefficients.  The famous  {\it Positivity Conjecture} 
says that  coefficients of the Laurent polynomials corresponding to  variables $x_i$ are always non-negative integers, see [Fomin \& Zelevinsky 2002]  \cite{FoZe1}.  
 A general form of the Positivity Conjecture  was proved by [Lee \& Schiffler 2015]   \cite{LeSch1}   using
 a  clever combinatorial formula for the variables $x_i$.

{\it Cluster $C^*$-algebras} ${\Bbb A}(\mathbf{x}, B)$  are a class of non-commutative rings introduced in 
\cite{Nik1}. The ${\Bbb A}(\mathbf{x}, B)$  is an {\it $AF$-algebra}
 given by the Bratteli diagram [Bratteli 1972]  \cite{Bra1}; 
 such a diagram is  obtained from a  mutation tree of the initial seed $(\mathbf{x}, B)$
 modulo an equivalence relation between the seeds  lying at  the same level,   see 
Section 2.2. 
(We refer the reader to Figures 1 and 2 for an immediate  example of such algebras.)

The aim of our note is the $K$-theory of   the $AF$-algebra  ${\Bbb A}(\mathbf{x}, B)$.
Namely,   the  {\it ordered} abelian group is a pair $(G,G^+)$ 
consisting of an abelian group $G$ and a semigroup $G^+\subset G$
of positive elements of $G$;  the order $\le$ on $G$ is defined by the 
positive cone $G^+$,  i.e. $a\le b$ if and only if $b-a\in G^+$.
An {\it order-isomorphism} $\cong$  between   $(G,G^+)$
and $(H,H^+)$ is an isomorphism $\varphi: G\to H$,  such that 
$\varphi(G^+)=H^+$.  Denote by   $K_0({\Bbb A}(\mathbf{x}, B))$
the $K_0$-group of the $AF$-algebra ${\Bbb A}(\mathbf{x}, B)$  
and by   $K_0^+({\Bbb A}(\mathbf{x}, B))\subset K_0({\Bbb A}(\mathbf{x}, B))$
its Grothendieck semigroup  [Effros 1981, Chapter 8]  \cite{E}.
The pair   $(K_0({\Bbb A}(\mathbf{x}, B)),  K_0^+({\Bbb A}(\mathbf{x}, B)))$
is an invariant of  {\it Morita equivalence}  of the $AF$-algebra ${\Bbb A}(\mathbf{x}, B)$
[Elliott 1976]   \cite{Ell1}.  In view of the Laurent phenomenon,  let  ${\cal A}_{{\bf add}}(\mathbf{x}, B)$ 
be an additive group of the cluster algebra  ${\cal A}(\mathbf{x}, B)$;  
let    ${\cal A}_{{\bf add}}^+(\mathbf{x}, B)$  be  a   
semigroup inside the ${\cal A}_{{\bf add}}(\mathbf{x}, B)$  consisting of  the Laurent  polynomials 
with {\it positive} coefficients.  The pair $({\cal A}_{{\bf add}}(\mathbf{x}, B), {\cal A}_{{\bf add}}^+(\mathbf{x}, B))$
is an abelian group with order.  The order $a>b$ is defined between two elements $a,b\in {\cal A}_{{\bf add}}(\mathbf{x}, B)$ if and only
if $a-b\in {\cal A}_{{\bf add}}^+(\mathbf{x}, B)$ .  Our main result   can be formulated as follows. 
\begin{thm}\label{thm1}
$(K_0({\Bbb A}(\mathbf{x}, B)), ~K_0^+({\Bbb A}(\mathbf{x}, B))
\cong ({\cal A}_{{\bf add}}(\mathbf{x}, B),  ~{\cal A}_{{\bf add}}^+(\mathbf{x}, B)).$
\end{thm}
An application of theorem \ref{thm1} is as follows.  Recall that the {\it dimension group}
is a triple $(G, G^+, \Gamma)$  consisting of an abelian group $G$,  a semigroup  
of positive elements  $G^+\subset G$ and a scale $\Gamma\subseteq  G^+$, i.e. a generating,
hereditary and directed subset of $G^+$    [Effros 1981, Chapter 7]  \cite{E}. 
For instance, the $\Gamma\cong G^+$ is a scale called {\it stable};  thus the  pair 
$(G,G^+)$ is a special case of the dimension group. 
An order-isomorphism $\cong$  between dimension groups    $(G,G^+, \Gamma)$
and $(H,H^+, \Gamma')$ is an isomorphism $\varphi: G\to H$,  such that 
$\varphi(G^+)=H^+$ and $\varphi(\Gamma)=\Gamma'$. 
Denote by $\Gamma\subset K_0^+({\Bbb A}(\mathbf{x}, B))$   the set of  the Murray-von Neumann  
equivalence classes  of projections in the algebra   ${\Bbb A}(\mathbf{x}, B)$.  
It is known, that the triple  $(K_0({\Bbb A}(\mathbf{x}, B)),  K_0^+({\Bbb A}(\mathbf{x}, B)), \Gamma)$
is an  invariant of  the {\it isomorphism}  class of the 
$AF$-algebra  ${\Bbb A}(\mathbf{x}, B)$  [Elliott 1976]   \cite{Ell1}.  
It is not hard to observe, that the set  $X=\{x_i\}_{i=1}^{\infty}$ of all variables $x_i$  in the cluster algebra
 ${\cal A}_{{\bf add}}(\mathbf{x}, B)$ is a scale,  since it is  a generating, hereditary and directed subset 
 of ${\cal A}^+_{{\bf add}}(\mathbf{x}, B)$. 
 Notice that   choosing a different initial seed $(\mathbf{x}, B)$  for the Laurent expansion 
of  variables $x_i$ yields a new scale $X'$,  such that 
$({\cal A}_{{\bf add}}(\mathbf{x}, B),  {\cal A}_{{\bf add}}^+(\mathbf{x}, B), X)$ 
 $\cong$  $({\cal A}_{{\bf add}}(\mathbf{x}, B),  {\cal A}_{{\bf add}}^+(\mathbf{x}, B), X')$. 
  But $X \subseteq  {\cal A}_{{\bf add}}^+(\mathbf{x}, B)$ for any 
dimension group;   therefore theorem \ref{thm1}  implies a new proof of  the Positivity Conjecture 
for the cluster algebras.
\begin{cor}\label{cor1}
The coefficients of the Laurent polynomials corresponding to the cluster variables $x_i$
are non-negative integers.   
\end{cor}
The article is organized as follows.  The preliminary  facts are 
introduced in Section 2.  Theorem \ref{thm1} and corollary \ref{cor1} are proved
in Section 3.  In Section 4 we consider an example of the cluster $C^*$-algebra 
${\Bbb A}(1,1)$ coming  from triangulation of an annulus with one marked point on each boundary  component.

\section{Preliminaries}
This section is a brief review of the  $AF$-algebras, cluster
$C^*$-algebras and Mundici dimension groups.  For a general review of $C^*$-algebras we refer the reader to 
[Murphy   1990]  \cite{M}.   The $AF$-algebras were introduced in   
 [Bratteli 1972]  \cite{Bra1}.  The general $K$-theory of $C^*$-algebras 
 is covered in [R\o rdam, Larsen \& Laustsen  2000]  \cite{RLL}  and $K$-theory of the 
 $AF$-algebras in  [Effros 1981] \cite{E}.  Cluster $C^*$-algebras 
 were the subject of \cite{Nik1}.  Mundici dimension groups were introduced 
 by  [Mundici 1988] \cite{Mun1}.

\subsection{$AF$-algebras and dimension groups}
A {\it $C^*$-algebra} is an algebra $A$ over $\mathbf{C}$ with a norm
$a\mapsto ||a||$ and an involution $a\mapsto a^*$ such that
it is complete with respect to the norm and $||ab||\le ||a||~ ||b||$
and $||a^*a||=||a||^2$ for all $a,b\in A$.
Any commutative $C^*$-algebra is  isomorphic
to the algebra $C_0(X)$ of continuous complex-valued
functions on some locally compact Hausdorff space $X$; 
otherwise, $A$ can be thought of as  a noncommutative  topological
space.

An {\it $AF$-algebra}  (Approximately Finite $C^*$-algebra) is defined to
be the  norm closure of a dimension-increasing  sequence of   finite dimensional
$C^*$-algebras $M_n$,  where  $M_n$ is the $C^*$-algebra of the $n\times n$ matrices
with entries in $\mathbf{C}$. Here the index $n=(n_1,\dots,n_k)$ represents
the  semi-simple matrix algebra $M_n=M_{n_1}\oplus\dots\oplus M_{n_k}$.
The ascending sequence mentioned above  can be written as 
\begin{equation}\label{eq2}
M_1\buildrel\rm\varphi_1\over\longrightarrow M_2
   \buildrel\rm\varphi_2\over\longrightarrow\dots,
\end{equation}
where $M_i$ are the finite dimensional $C^*$-algebras and
$\varphi_i$ the homomorphisms between such algebras.  
The homomorphisms $\varphi_i$ can be arranged into  a graph as follows. 
Let  $M_i=M_{i_1}\oplus\dots\oplus M_{i_k}$ and 
$M_{i'}=M_{i_1'}\oplus\dots\oplus M_{i_k'}$ be 
the semi-simple $C^*$-algebras and $\varphi_i: M_i\to M_{i'}$ the  homomorphism. 
(To keep it simple, one can assume that $i'=i+1$.) 
One has  two sets of vertices $V_{i_1},\dots, V_{i_k}$ and $V_{i_1'},\dots, V_{i_k'}$
joined by  $b_{rs}$ edges  whenever the summand $M_{i_r'}$ contains $b_{rs}$
copies of the summand $M_{i_s}$ under the embedding $\varphi_i$. 
As $i$ varies, one obtains an infinite graph called the  {\it Bratteli diagram} of the
$AF$-algebra.  The matrix $B=(b_{rs})$ is known as  a {\it partial multiplicity} matrix;
an infinite sequence of $B_i$ defines a unique $AF$-algebra.

For a unital $C^*$-algebra $A$, let $V(A)$
be the union (over $n$) of projections in the $n\times n$
matrix $C^*$-algebra with entries in $A$;
projections $p,q\in V(A)$ are {\it equivalent} if there exists a partial
isometry $u$ such that $p=u^*u$ and $q=uu^*$. The equivalence
class of projection $p$ is denoted by $[p]$;
the equivalence classes of orthogonal projections can be made to
a semigroup by putting $[p]+[q]=[p+q]$. The Grothendieck
completion of this semigroup to an abelian group is called
the  $K_0$-group of the algebra $A$.
The functor $A\to K_0(A)$ maps the category of unital
$C^*$-algebras into the category of abelian groups, so that
projections in the algebra $A$ correspond to a positive
cone  $K_0^+\subset K_0(A)$ and the unit element $1\in A$
corresponds to an order unit $u\in K_0(A)$.
The ordered abelian group $(K_0,K_0^+,u)$ with an order
unit  is called a {\it dimension group};  an order-isomorphism
class of the latter we denote by $(G,G^+)$.

If ${\Bbb A}$ is an $AF$-algebra, then its dimension group
$(K_0({\Bbb A}), K_0^+({\Bbb A}), u)$ is a complete isomorphism
invariant of algebra ${\Bbb A}$ [Elliott 1976]   \cite{Ell1}.  
The order-isomorphism  class $(K_0({\Bbb A}), K_0^+({\Bbb A}))$
 is an invariant of  {\it Morita equivalence} of the  algebra 
${\Bbb A}$,  i.e.  an isomorphism class in the category of 
finitely generated projective modules over ${\Bbb A}$.

The {\it scale} $\Gamma$ is a subset of  $K_0^+({\Bbb A})$
which is generating, hereditary and directed, i.e. 
(i)  for each $a\in K_0^+({\Bbb A})$ there exist $a_1,\dots,a_r\in\Gamma({\Bbb A})$,
such that $a=a_1+\dots+a_r$;
(ii)  if $0\le a\le b\in \Gamma$, then $a\in\Gamma$;
(iii) given $a,b\in\Gamma$ there exists $c\in\Gamma$, such that
$a,b\le c$.  If $u$ is an order unit,  then the set $\Gamma:=\{a\in K_0^+({\Bbb A})~|~ 0\le a\le u\}$
is a scale;  thus the dimension group of algebra ${\Bbb A}$ can be written 
in the form $(K_0({\Bbb A}), K_0^+({\Bbb A}), \Gamma)$.

\subsection{Cluster $C^*$-algebras}
Let $T_n$ be an oriented tree whose vertices correspond to the seeds 
$(\mathbf{x}, B)$ and outgoing edges correspond to mutations in 
direction $1\le k\le n$.   Notice that the  tree $T_n$  of 
a cluster algebra ${\cal A}(\mathbf{x}, B)$ has a grading by  levels,
i.e. the minimal distance from the root of $T_n$. 
We shall say that a pair of  clusters $\mathbf{x}$ and $\mathbf{x}'$ with exchange matrices
$B$ and $B'$ are   {\it $\ell$-equivalent},    if:

\medskip
  (i)  $\mathbf{x}$ and $\mathbf{x}'$ lie at the same level;   

\smallskip
(ii) $\mathbf{x}$ and $\mathbf{x}'$  coincide
modulo a cyclic permutation of variables $x_i$;

\smallskip
(iii)  $B=B'$.  

\medskip\noindent
It is not hard to see that $\ell$ is an equivalence relation on the set 
of vertices of graph $T_n$. 
\begin{dfn}\label{cluster}
By  a cluster $C^*$-algebra ${\Bbb A}(\mathbf{x}, B)$ 
 one understands an $AF$-algebra given by
the Bratteli diagram ${\goth B}(\mathbf{x}, B)$ of the form:
\begin{equation}
{\goth B}(\mathbf{x}, B):= T_n ~\hbox{\bf mod} ~\ell. 
\end{equation}
\end{dfn}
\begin{rmk}
\textnormal{
Notice  that  the graph  ${\goth B}(\mathbf{x}, B)$ is no longer a
tree;  the cycles of ${\goth B}(\mathbf{x}, B)$ appear after 
gluing together vertices lying at the same level  of the tree
according to the relation $\ell$.   
The ${\goth B}(\mathbf{x}, B)$ is not a regular graph, since the valency
of vertices can vary.  
However,  the ${\goth B}(\mathbf{x}, B)$ is always a Bratteli diagram, since 
it is obtained from a regular tree by an addition of extra edges and subsequent 
contraction of the respective edges. 
Notice also that 
 the ${\goth B}(\mathbf{x}, B)$  is a finite graph if and only if 
${\cal A}(\mathbf{x}, B)$ is a  finite  cluster algebra. 
}
\end{rmk}
\begin{exm}
\textnormal{
Let $\mathbf{x}=(x_1,x_2)$ and
\begin{equation}
B=\left(\matrix{0 & 2\cr -2 & 0}\right).  
 \end{equation}
The cluster algebra ${\cal  A}(\mathbf{x}, B)$  is associated to an ideal  triangulation 
 of an annulus  with one marked point  on each boundary component,
see [Fomin,  Shapiro  \& Thurston  2008, Example 4.4]  \cite{FoShaThu1}.  
The exchange relations (\ref{eq1})  in this case can be written as $x_{i-1}x_{i+1}=1+x_i^2$
and  $B'=-B$.  It is easy  to verify using 
definition \ref{cluster},  that   the Bratteli diagram $T_2 ~\hbox{\bf mod}~\ell$  of the corresponding 
cluster $C^*$-algebra  ${\Bbb A}(\mathbf{x}, B)$ is given by 
Figure 3.  We refer the reader to Section 4  for an extended  discussion of the properties 
of such an algebra.  
}
\end{exm}
\begin{exm}\label{exm1}
{\bf (\cite{Nik1})}
\textnormal{
Let $\mathbf{x}=(x_1,x_2,x_3)$ and
\begin{equation}
 B=\left(\matrix{0 & 2 &-2\cr -2 & 0 & 2\cr 2 & -2 & 0}\right).   
\end{equation}
The cluster algebra ${\cal  A}(\mathbf{x}, B)$ is called {\it Markov's};
it is associated to an ideal  triangulation of the hyperbolic torus with a cusp,
see e.g. [Williams 2014] \cite{Wil1}.   The Bratteli diagram $T_3 ~\hbox{\bf mod} ~\ell$
of the cluster $C^*$-algebra  ${\Bbb A}(\mathbf{x}, B)$ is shown in 
Figure 1. 
(The corresponding mutation tree $T_3$ and the equivalence  classes of relation $\ell$ 
are given in full detail in \cite{Nik1},  Figure 4.) 
 The algebra  ${\Bbb A}(\mathbf{x}, B)$ has a  non-trivial  primitive spectrum being  
 isomorphic to an  $AF$-algebra ${\goth M}_1$ introduced 
by  [Mundici 1988]  \cite{Mun1};  for a general theory we refer the reader 
to the monograph [Mundici 2011]   \cite{MU}. 
}
\end{exm}
\begin{figure}
\begin{picture}(100,100)(-160,130)

\put(50,200){\circle{3}}
\put(50,179){\circle{3}}
\put(20,179){\circle{3}}
\put(80,179){\circle{3}}


\put(3,166){\circle{3}}
\put(20,166){\circle{3}}
\put(36,166){\circle{3}}
\put(50,166){\circle{3}}
\put(65,166){\circle{3}}

\put(80,166){\circle{3}}
\put(98,166){\circle{3}}


\put(49,199){\vector(-3,-2){29}}
\put(51,199){\vector(3,-2){29}}
\put(50,200){\vector(0,-1){20}}


\put(19,178){\vector(-3,-2){15}}
\put(21,178){\vector(3,-2){15}}
\put(20,178){\vector(0,-1){10}}


\put(50,178){\vector(-3,-2){15}}
\put(50,178){\vector(3,-2){15}}
\put(50,178){\vector(0,-1){10}}


\put(79,178){\vector(-3,-2){15}}
\put(81,178){\vector(3,-2){15}}
\put(80,178){\vector(0,-1){10}}


\put(1.5,165){\vector(-1,-1){13}}
\put(4,164){\vector(2,-3){7}}
\put(2.5,164){\vector(0,-1){10}}

\put(19,164){\vector(-1,-1){10}}
\put(21,164){\vector(1,-1){10}}
\put(20,164){\vector(0,-1){10}}

\put(49,164){\vector(-1,-1){10}}
\put(51,164){\vector(1,-1){10}}
\put(50,164){\vector(0,-1){10}}

\put(79,164){\vector(-1,-1){10}}
\put(81,164){\vector(1,-1){10}}
\put(80,164){\vector(0,-1){10}}

\put(98,164){\vector(-1,-1){10}}
\put(100,164){\vector(1,-1){10}}
\put(98,164){\vector(0,-1){10}}


\put(65,164){\vector(-1,-1){10}}
\put(65,164){\vector(1,-1){10}}
\put(65,164){\vector(0,-1){10}}

\put(36,164){\vector(-1,-1){10}}
\put(36,164){\vector(1,-1){10}}
\put(36,164){\vector(0,-1){10}}


\end{picture}
\caption{The Bratteli diagram of  Markov's  cluster $C^*$-algebra.} 
\end{figure}
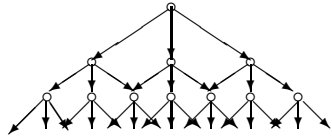

\subsection{Mundici dimension groups}
A broad class of dimension groups has been introduced by [Mundici 1988] \cite{Mun1}. 
We shall use such groups as the main  technical tool in proof of theorem \ref{thm1}. 
We refer the reader to   [Mundici 1988] \cite{Mun1} and [Mundici 2011]   \cite{MU}
for a detailed account.  

A {\it lattice-ordered} (or $\ell$-group) is a structure $(G,+,-,0,\vee,\wedge)$ such that 
$(G,+,-,0)$ is an abelian group, $(G,\vee,\wedge)$ is a lattice, and $x+(y\vee z)=(x+y)\vee (x+z)$
for all $x,y,z\in G$. An {\it order unit} in a partially ordered group $G$ is an element $u\ge 0$ such
that for each $x\in G$ there is an integer $n\ge 0$ with $x<nu$.   
A {\it unital $\ell$-group} is an $\ell$-group with distinguished order unit.  

The function $f: [0, 1]^n\to {\Bbb R}$  is called  a {\it McNaughton
function} over $[0, 1]^n$  iff $f$  is continuous and there are a finite number of linear
functions: 
\begin{equation}\label{equ24}
\left\{
\begin{array}{lll}
\alpha_1 &=&  b_1 +a_{11}x_1+\dots+a_{1n}x_n\\
\alpha_2 &=&  b_2 +a_{21}x_1+\dots+a_{2n}x_n\\
\vdots\\
\alpha_m &=&  b_1 +a_{m1}x_1+\dots+a_{mn}x_n,
\end{array}
\right.
\end{equation}
where all $a_{ij}$  and $b_i$  are integers, such that for 
every $(x_1, \dots,  x_n)\in [0, 1]^n$  there is $i\in \{1, \dots, m\}$
with $f(x_1, \dots, x_n) = \alpha_i(x_1,\dots,x_n)$, see  
 [Mundici 1988] \cite{Mun1} and [Mundici 2011]   \cite{MU}. 
In other words, the McNaughton function is a piecewise linear 
function with integer coefficients.  It is easy to see,  that  the set of all McNaughton
functions over $[0,1]$  is an $\ell$-group with the pointwise operations $+,-,\vee,\wedge$ of ${\Bbb R}$
and with the constant function $1$ as the distinguished order unit.
The {\it Mundici dimension group} ${\cal M}_n$ is an $\ell$-group defined by  the McNaughton functions over $[0,1]^n$. 
\begin{thm}\label{thm2}
{\bf (\cite{Mun1}, \cite{MU})}
$(K_0({\Bbb A}(\mathbf{x}, B)), K_0^+({\Bbb A}(\mathbf{x}, B)),u)\cong ({\cal M}_n,1)$,
where ${\cal M}_n$ is defined by  a subset of  all  McNaughton functions over $[0,1]^n$.   
\end{thm}
\begin{rmk}
\textnormal{
Theorem \ref{thm2} for $n=1$ was proved in  [Mundici 1988] \cite{Mun1}.
In particular, the Markov cluster $C^*$-algebra ${\Bbb A}(\mathbf{x}, B)$ in Figure 1 has 
the dimension group ${\cal M}_1$. 
By an extension of the argument of  [Mundici 2011]   \cite{MU},  
one can prove Theorem \ref{thm2} for $n\ge 1$. 
}
\end{rmk}

\section{Proofs}
\subsection{Proof of theorem \ref{thm1}}
We shall split the proof into a series of lemmas.
\begin{lem}\label{lem0}
The ordered abelian group  $({\cal A}_{{\bf add}}(\mathbf{x}, B),  ~{\cal A}_{{\bf add}}^+(\mathbf{x}, B))$
is a dimension group  with the  stable scale $\Gamma\cong {\cal A}_{{\bf add}}^+(\mathbf{x}, B)$. 
\end{lem}
{\it Proof.} 
Recall that an ordered abelian group $(G, G^+)$ satisfies the {\it Riesz interpolation property}, 
if given $\{a_i, b_j\in G ~|~ a_i\le b_j ~\hbox{for} ~i,j=1,2\}$  there exists $c\in G$,  such that
\begin{equation}\label{riesz}
a_i\le c\le b_j. 
\end{equation}
Let us show that the ordered group  $({\cal A}_{{\bf add}}(\mathbf{x}, B),  ~{\cal A}_{{\bf add}}^+(\mathbf{x}, B))$
satisfies the Riesz interpolation property.  Indeed,  if $a_i=\sum A_i\mathbf{x}^i$ and 
$b_j=\sum B_i\mathbf{x}^i$ are the Laurent polynomials of $a_i,b_j\in {\cal A}_{{\bf add}}(\mathbf{x}, B)$,
then one can choose $c=\sum C_i\mathbf{x}^i$ such that $C_i=A_i$ if $A_i\ne 0$ and $0< C_i<B_i$
if $A_i=0$.  Clearly,  the condition (\ref{riesz}) is satisfied. 

By the Effros-Handelman-Shen Theorem,  a countable ordered abelian group is  a dimension group
if and only if it satisfies the Riesz interpolation property,  see  [Effros 1981, Theorem 3.1] \cite{E}. 
Thus $({\cal A}_{{\bf add}}(\mathbf{x}, B),  ~{\cal A}_{{\bf add}}^+(\mathbf{x}, B))$ is a dimension 
group with the stable scale. Lemma \ref{lem0}  is proved.

\begin{lem}\label{lem3}
The exists a canonical isomorphism $\varphi$ between the abelian groups
$K_0({\Bbb A}(\mathbf{x}, B))$ and ${\cal A}_{{\bf add}}(\mathbf{x}, B)$.
 \end{lem}
{\it Proof.} 
The idea is to construct an isomorphism $\varphi: {\cal M}_n\to {\cal A}_{{\bf add}}(\mathbf{x}, B)$,
where ${\cal M}_n$ is the Mundici dimension group. The rest of the proof will follow from  Theorem \ref{thm2}.  

We assume without loss of generality, that the linear functions $\alpha_i$ in the set (\ref{equ24})
constitute the Schauder-type  basis for $[0,1]^n$.
\footnote{Such a choice of $\alpha_i$ provides injectivity of our construction. I am grateful to the referee for pointing out this fact to me.}
 We shall assign to each  $\alpha_i =  b_i +a_{i1}x_1+\dots+a_{in}x_n$
a Laurent monomial, which is  a generator of  the group   ${\cal A}_{{\bf add}}(\mathbf{x}, B)$. Roughly speaking, this can be done
by an ``exponentiation'' of the variables $x_i$.  

Indeed, consider a map $\varphi$ acting by the formula:
\begin{equation}\label{equ8}
b_i +a_{i1}x_1+a_{i2}x_2+\dots+a_{in}x_n\mapsto b_i x_1^{a_{i1}} x_2^{a_{i2}}\dots  x_n^{a_{in}}, 
\quad 1\le i\le m, 
\end{equation}
where $a_{ij}\in \mathbf{Z}$ and  $b_i\in \mathbf{Z}$.  

The map $\varphi$ sends the McNaughton function $f:[0,1]^n\to {\Bbb R}$ to a Laurent
polynomial according to the formula:
\begin{equation}
f\mapsto \sum_{i=1}^m b_i x_1^{a_{i1}} x_2^{a_{i2}}\dots  x_n^{a_{in}}\in \mathbf{Z}[\mathbf{x}^{\pm 1}].
\end{equation}

It is verified directly, that  pointwise addition of the McNaughton functions maps to an addition of the Laurent polynomials.
One can see that our construction provides an injective mapping into the ring of Laurent polynomials.  
Working backwards our construction, it can be proved that the mapping is also  surjective.   
 Thus one gets an isomorphism of the abelian groups:
\begin{equation}
\varphi: {\cal M}_n\to {\cal A}_{{\bf add}}(\mathbf{x}, B).
\end{equation}

On the other hand, it follows from Theorem \ref{thm2} that  ${\cal M}_n\cong K_0({\Bbb A}(\mathbf{x}, B))$,
where the set of the   McNaughton functions over $[0,1]^n$ is defined by the algebra ${\Bbb A}(\mathbf{x}, B)$. 
Thus one obtains an isomorphism of the abelian groups:
\begin{equation}
\varphi: K_0({\Bbb A}(\mathbf{x}, B)) \to {\cal A}_{{\bf add}}(\mathbf{x}, B).
\end{equation}

Lemma \ref{lem3} is proved.

\begin{rmk}
\textnormal{
Using the McNaughton functions over $[0,1]^n$, one can  see that for the 
finite cluster algebras the group $K_0({\Bbb A}(\mathbf{x}, B))$ is isomorphic to 
 a direct sum of finitely many copies of $\mathbf{Z}$.
This fact  implies that the corresponding tracial simplex is spanned by $n$ extremal traces.
}
\end{rmk}
\begin{rmk}
\textnormal{
Lemma \ref{lem3} implies  that the group $K_0({\Bbb A}(\mathbf{x}, B))$ has the 
natural structure of a commutative ring, since $K_0({\Bbb A}(\mathbf{x}, B))\subset \mathbf{Z}[\mathbf{x}^{\pm 1}]$. 
It is an interesting question to find an interpretation of the product in terms of the $K$-theory.  
 }
\end{rmk}
\begin{lem}\label{lem4}
The isomorphism $\varphi$ is order-preserving, i.e. 
$$(K_0({\Bbb A}(\mathbf{x}, B)), ~K_0^+({\Bbb A}(\mathbf{x}, B))
\cong ({\cal A}_{{\bf add}}(\mathbf{x}, B),  ~{\cal A}_{{\bf add}}^+(\mathbf{x}, B)).$$
\end{lem}
{\it Proof.} 
In view of Theorem \ref{thm2}, it is sufficient to show that
\begin{equation}
({\cal M}_n, 1)\cong ({\cal A}_{{\bf add}}(\mathbf{x}, B),  ~{\cal A}_{{\bf add}}^+(\mathbf{x}, B))
\end{equation}
are isomorphic dimension groups. 

The semi-group ${\cal M}_n^+$ of positive elements of the Mundici dimension group 
$({\cal M}_n, 1)$ consists of all piecewise linear functions with $b_i>0$. 
Likewise, the semigroup ${\cal A}_{{\bf add}}^+(\mathbf{x}, B)$ consists of the Laurent polynomials 
with $b_i>0$. 

On the other hand, formula (\ref{equ8}) says that $\varphi$ sends the coefficient $b_i$ 
into the coefficient $b_i$ of the Laurent monomial. 
Thus one gets the equality:
\begin{equation}\label{equ13}
\varphi({\cal M}_n^+)={\cal A}_{{\bf add}}^+(\mathbf{x}, B).
\end{equation}
In other words, the isomorphism $\varphi: {\cal M}_m\to {\cal A}_{{\bf add}}(\mathbf{x}, B)$ preserves  
the semi-group of positive elements of the respective dimension groups.

Lemma \ref{lem4} follows from (\ref{equ13}) and  Theorem \ref{thm2}.

\bigskip
Theorem \ref{thm1} follows from lemma \ref{lem4}.

\bigskip
\begin{rmk}
{\bf (\cite{Nik1})}
\textnormal{
Theorem \ref{thm1}  implies that   the category of cluster algebras 
can be embedded into the category of dimension groups $(G,G^+)$ with the stable scale.
The following (partial)  characterization of  cluster algebras in terms of the dimension 
groups is true:   The cluster algebras correspond to the dimension groups with
a non-trivial spectrum $Prim~(G, G^+) \cong \{\mathbf{R}^n ~|~n\ge 1\}$,
where   $Prim~(G, G^+)$ is  the space of primitive ideals of  $(G, G^+)$
endowed with the Jacobson topology. 
}
\end{rmk}

\subsection{Proof of corollary \ref{cor1}}
Let $\psi$ be an inverse of the map $\varphi$ constructed in lemma \ref{lem3}.
We shall fix an isomorphism class of the $AF$-algebra   ${\Bbb A}(\mathbf{x}, B)$
and consider the corresponding dimension group 
$(K_0({\Bbb A}(\mathbf{x}, B)),  K_0^+({\Bbb A}(\mathbf{x}, B)), \Gamma)$. 
In view of theorem \ref{thm1}, we have:
\begin{equation}
\left\{
\begin{array}{lll}
{\cal A}_{{\bf add}}(\mathbf{x}, B) &=&\psi(K_0({\Bbb A}(\mathbf{x}, B)))\\
{\cal A}^+_{{\bf add}}(\mathbf{x}, B) &=&\psi(K_0^+({\Bbb A}(\mathbf{x}, B))). 
\end{array}
\right.
\end{equation}
Since $\Gamma\subseteq K_0^+({\Bbb A}(\mathbf{x}, B))$,  
one gets a scale $\psi(\Gamma)\subseteq  {\cal A}^+_{{\bf add}}(\mathbf{x}, B)$
in the cluster algebra ${\cal A}_{{\bf add}}(\mathbf{x}, B)$.  

On the other hand, it is verified directly that the set $X=\{x_i\}_{i=1}^{\infty}$ of all 
cluster variables $x_i$ is a scale,  since it is a  generating,  hereditary and directed subset 
 of ${\cal A}_{{\bf add}}(\mathbf{x}, B)$.  But given isomorphism class of 
 algebra ${\Bbb A}(\mathbf{x}, B)$ can define  only one scale on the cluster
 algebra  ${\cal A}_{{\bf add}}(\mathbf{x}, B)$;  thus $X\cong  \psi(\Gamma)$.
 It remains to recall  that  $\psi(\Gamma) \subseteq  {\cal A}^+_{{\bf add}}(\mathbf{x}, B)$
 and therefore $X \subseteq  {\cal A}^+_{{\bf add}}(\mathbf{x}, B)$.
 In other words,  the coefficients of the Laurent polynomials corresponding to the cluster variables $x_i$
are non-negative integers. 
 Corollary \ref{cor1} is proved.

\section{An example}
To illustrate theorem \ref{thm1}, we shall consider a cluster $C^*$-algebra ${\Bbb A}(1,1)$ 
associated to a triangulation of an annulus  with one marked point
on each boundary component,
see [Fomin,  Shapiro  \& Thurston  2008, Example 4.4]  \cite{FoShaThu1};
we shall keep the original notation of cited paper.

\subsection{Cluster $C^*$-algebra ${\Bbb A}(1,1)$}
Let   ${\goth A}=\{z=x+iy\in\mathbf{C} ~|~ r\le |z|\le R\}$ be an annulus 
in the complex plane such that  $r<R$.  
Recall that  the Riemann surfaces  ${\goth A}$ and ${\goth A}'$  are  isomorphic if and only if     
$R/r=R'/r'$;  the real number  $t=R/r$ is called a {\it modulus} of ${\goth A}$.
By $T_{\goth A}=\{t\in \mathbf{R} ~|~ t>1\}$  we understand the Teichm\"uller
space of ${\goth A}$,  i.e. the space of all  Riemann surfaces ${\goth A}$ endowed with a 
natural  topology.   The cluster algebra ${\cal A}(\mathbf{x}, B_T)$ of rank two given by a
matrix: 
\begin{equation}\label{eq15}
B_T=\left(\matrix{0 & 2\cr -2 & 0}\right)  
\end{equation}
is the coordinate ring of  $T_{\goth A}$ [Fomin,  Shapiro  \& Thurston  2008, Example 4.4]  \cite{FoShaThu1}; 
 the  ${\cal A}(\mathbf{x}, B_T)$ is related to the Penner coordinates on the space  $T_{\goth A}$ 
 corresponding  to an ideal  triangulation   $T$ of  ${\goth A}$  with one marked point  on each boundary component of  ${\goth A}$
  [Williams 2014, Section 3] \cite{Wil1}.

By  ${\Bbb A}(1,1):={\Bbb A}(\mathbf{x}, B_T)$ we shall understand a cluster
$C^*$-algebra given by matrix $B_T$;   the reader is encouraged to verify 
using formulas (\ref{eq1})  that the Bratteli diagram of  ${\Bbb A}(1,1)$
has the form of a  Pascal triangle shown in  Figure 3. 
(The ${\Bbb A}(1,1)$ is  the so-called 
{\it GICAR  algebra}  [Effros 1980, p.13(e)]  \cite{E};  such an algebra has a rich set  of 
ideals [Bratteli 1972, Section 5.5]  \cite{Bra1}.)

\begin{figure}
\begin{picture}(100,100)(-150,150)

\put(50,200){\circle{3}}

\put(33,188){\circle{3}}
\put(67,188){\circle{3}}

\put(50,177){\circle{3}}
\put(16,177){\circle{3}}
\put(84,177){\circle{3}}

\put(-1,164){\circle{3}}
\put(34,164){\circle{3}}
\put(68,164){\circle{3}}
\put(103,164){\circle{3}}


\put(49,199){\vector(-3,-2){15}}
\put(51,199){\vector(3,-2){15}}

\put(32,187){\vector(-3,-2){15}}
\put(34,187){\vector(3,-2){15}}

\put(66,187){\vector(-3,-2){15}}
\put(68,187){\vector(3,-2){15}}


\put(14,175){\vector(-3,-2){15}}
\put(17,175){\vector(3,-2){15}}

\put(49,175){\vector(-3,-2){15}}
\put(51,175){\vector(3,-2){15}}

\put(83,175){\vector(-3,-2){15}}
\put(86,175){\vector(3,-2){15}}

\put(-10,155){$\dots$}
\put(27,155){$\dots$}
\put(64,155){$\dots$}
\put(101,155){$\dots$}

\end{picture}
\caption{Bratteli diagram of algebra ${\Bbb A}(1,1)$.} 
\end{figure}
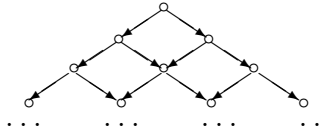

\bigskip\noindent
By $\{\sigma_t:  {\Bbb A}(1,1)\to {\Bbb A}(1,1) ~|~ t\in \mathbf{R}\}$ we denote a group of modular 
automorphisms constructed in  \cite{Nik1}, Section 4;  the $\sigma_t$  is 
generated by the  geodesic flow $T^t$ on  the space  $T_{\goth A}$.

The ${\Bbb A}(1,1)$  embeds into an UHF-algebra: 
\begin{equation}
M_{2^{\infty}}:=\bigotimes_{i=1}^{\infty} M_2(\mathbf{C}). 
\end{equation}
The $M_{2^{\infty}}$ is known as a  {\it CAR algebra};  
unlike the  ${\Bbb A}(1,1)$,  it  is a simple $AF$-algebra with  
the Bratteli diagram shown in  [Effros 1980, p.13(c1)]  \cite{E}.
The Powers product $\{\otimes_{i=1}^{\infty}\exp (\sqrt{-1}\left(\small\matrix{1 & 0\cr 0 & \lambda}\right))
~|~ 0<\lambda < 1\}$  defines a group of modular automorphisms 
$\{\sigma^t:   M_{2^{\infty}}\to M_{2^{\infty}} ~|~ t\in \mathbf{R}\}$;  
it is not hard to observe,  that  $\sigma^t \equiv \sigma_t$  on 
${\Bbb A}(1,1)$.

Recall that if $e_{ij}$ are the matrix units in $M_2(\mathbf{C})$, one can define a projection 
$e\in M_2(\mathbf{C})\otimes M_2(\mathbf{C})$ by the formula 
$e={1\over 1+t} (e_{11}\otimes e_{11}+te_{22}\otimes e_{22}+\sqrt{t}
(e_{12}\otimes e_{21}+e_{21}\otimes e_{12}))$,
where $t\in\mathbf{R}$ is a parameter. 
If  $\theta$ is the shift automorphism of the $UHF$-algebra $M_{2^{\infty}}$, 
then projections
$e_i:=\theta^i(e)\in M_{2^i}$ satisfy the following relations:
\begin{equation}\label{eq18}
\left\{
\begin{array}{ccc}
e_i e_j &=& e_j e_i, \quad\hbox{if} \quad |i-j|\ge 2\\
e_i e_{i\pm 1} e_i &=&  {t\over (1+t)^2} e_i
\end{array}
\right.
\end{equation}
 and  the Powers state $\varphi_t : M_{2^{\infty}}\to \mathbf{C}$ 
 satisfies the Jones equality:
\begin{equation}\label{eq19}
\varphi_t (we_{n+1})={t\over (1+t)^2}  ~\varphi_t(w), \quad\forall w\in M_{2^{n+1}},
\end{equation}
see  [Jones 1991, Section 5.6]  \cite{J1} for the details.   The $e_i$ 
generate the algebra $M_{2^{\infty}}$  and taking new generators $s_i$
such that $\sigma^t(s_i)=te_i-(1-e_i)$ one gets  a representation 
of the braid group $B_n=\{s_1,\dots, s_n ~|~ s_is_{i+1}s_i=s_{i+1}s_is_{i+1},
~s_is_j=s_js_i ~\hbox{if} ~|i-j|\ge 2\}$  in the algebra $M_{2^{n}}$.

\subsection{Jones Index Theorem}
As an application of theorem \ref{thm1}, one gets  a  short proof 
of the Jones Index Theorem in terms of the cluster algebras. 
\begin{cor}\label{corr2}
 Relations (\ref{eq18}) define a $C^*$-algebra if and only if the values of index ${(1+t)^2\over t}$ 
 belong to the set: 
\begin{equation}\label{eq33}
 [4, \infty) ~\bigcup  ~\{ 4\cos^2\left({\pi\over n}\right) ~|~n\ge 3\}.   
\end{equation}
\end{cor}
{\it Proof.} 
To find   admissible values of  parameter $t$, 
we shall use  a simple  analysis  of the cluster algebra  ${\cal A}(1,1)\cong K_0({\Bbb A}(1,1))$.  
Recall  that  algebra  ${\cal A}(1,1)$ has a unique canonical basis ${\cal B}$
consisting of the positive elements of  ${\cal A}(1,1)$,  i.e. the Laurent polynomials with positive 
integer coefficients;  the elements of $\cal B$ generate the whole algebra ${\cal A}(1,1)$. 
An explicit construction of $\cal B$ was given by  [Sherman \& Zelevinsky 2004, Theorem 2.8]  \cite{SheZe1}.  
Namely,
\begin{equation}\label{eq20}
{\cal B}=\{x_i^px_{i+1}^q ~|~ p,q\ge 0\} ~\bigcup ~\{T_n(x_1x_4-x_2x_3) ~|~ n\ge 3\},
\end{equation}
where $T_n(x)$ are the Chebyshev polynomials of the first kind.  
Since
\begin{equation}\label{eq21}
T_0=1 \quad\hbox{and}  \quad  T_n \left[{1\over 2}(t+t^{-1})\right]={1\over 2}(t^n+t^{-n}), 
\end{equation}
we shall look for a modulus $t$ such that  ${1\over 2}(t+t^{-1})=x_1x_4-x_2x_3$.
This is always possible since the Penner coordinates on the Teichm\"uller space $T_{\goth A}$
are given by the cluster $(x_1,x_2)$,  where each $x_i$   is a function of  modulus $t$  
[Williams 2014, Section 3] \cite{Wil1}.

\medskip
(i)  Since $t>1$,  it is easy to see by direct substitution that the values of index 
belong to the  interval:
\begin{equation}\label{eq27}
 (4, \infty).  
\end{equation}

\medskip
(ii)  To get discrete values,   we shall assume that  ${\cal A}(1,1)$
is  a finite cluster algebra,  i.e.  the number of  $x_i$ is finite.  
It is immediate that   $|{\cal B}|<\infty$  
and from the second series in  (\ref{eq20}) one obtains:
\begin{equation}\label{eq29}
T_n(x_1x_4-x_2x_3)=T_0=1
\end{equation}
for some integer $n\ge 1$.   But  $x_1x_4-x_2x_3= {1\over 2}(t+t^{-1})$
and using formula (\ref{eq21})  for the Chebyshev polynomials, one gets
an equation:
\begin{equation}\label{eq30}
t^n+t^{-n}=2
\end{equation}
for (possibly complex)  values of modulus $t$.  Since  (\ref{eq30}) 
is equivalent to the equation $t^{2n}-2t^n+1=(t^n-1)^2=0$,  one gets the 
$n$-th root of unity:
\begin{equation}\label{eq31}
t\in \{e^{2\pi i\over n} ~|~ n\ge 1\}. 
\end{equation}
However,  the index
\begin{equation}\label{eq32}
{(1+t)^2\over t}=t^{-1}+2+t=2\left[\cos~\left({2\pi\over n}\right)+1\right]=
4\cos^2\left({\pi\over n}\right)
\end{equation}
is a real number.   Thus  relations (\ref{eq18}) define a $C^*$-algebra. (We must exclude the 
case $n= 2$ corresponding to $t=-1$,  because otherwise one gets a division by zero in (\ref{eq18}).)

\medskip 
Bringing together   (\ref{eq27}) and (\ref{eq32}) one gets the conclusion 
of corollary \ref{corr2}. 
\begin{rmk}\label{rmk2}
\textnormal{
The finite  cluster algebras   corresponding to  the 
discrete moduli  come from a triangulation of the $n$-gons or the $n$-gons with one puncture,
see [Fomin,  Shapiro  \& Thurston  2008, Table 1]  \cite{FoShaThu1};
such  algebras are classified by their   Coxeter-Dynkin diagrams of type 
$A_{n-3}$ and $D_n$, respectively.  As explained,  the 
${\Bbb  A}(1,1)$  is a  finite-dimensional  $C^*$-algebra having   the  Bratteli 
diagram   similar to one  shown in   [Jones 1991, pp. 37-38]  \cite{J1}.  
}
\end{rmk}

\subsection{Dimension group of the GICAR algebra}
We shall use theorem \ref{thm1} to calculate  a dimension group of the algebra 
${\Bbb A}(1,1)$. 
\begin{cor}\label{cor2}
 $(K_0({\Bbb A}(1,1)), K_0^+({\Bbb A}(1,1)), u)\cong (\mathbf{Z}[x], P^+(0,1), u),$
where $P^+(0,1)$ is the semigroup of all positive-definite polynomials 
on the interval $(0,1)$. 
\end{cor}
{\it Proof.}  It is known that the Chebyshev polynomials of the first kind $T_n(x)$ lie in
a basis ${\cal B}$ of the cluster algebra  ${\cal A}(1,1)$  [Sherman \& Zelevinsky 2004, Theorem 2.8]  \cite{SheZe1}. 
For each $0\le k\le n$,  we shall introduce a new basis ${\cal B}'$ in    ${\cal A}(1,1)$
comprising the elements:
\begin{equation}
T_1^k(x)(T_0(x)-T_1(x))^{n-k}=x^k(1-x)^{n-k}.
\end{equation}
On the other hand, the Bratteli diagram in Figure 3 says that 
the group $K_0({\Bbb A}(1,1))$ is generated by the (equivalence classes of)
projections $[e_k^n]$ subject to the relations:
\begin{equation}\label{eq36}
[e_k^n]=[e_k^{n+1}]+[e_{k+1}^{n+1}].
\end{equation}
Take a representation $\rho$ of $K_0({\Bbb A}(1,1))$ in the 
cluster algebra  ${\cal A}(1,1)$ given by the formula:
\begin{equation}
\rho([e_k^n])=x^k(1-x)^{n-k}, \quad 0\le k\le n.
\end{equation}
The reader can verify that relations (\ref{eq36}) are satisfied.   It is easy to see, that $x^k(1-x)^{n-k}$ 
are generators of the polynomial ring $\mathbf{Z}[x]$  and the rest of the 
proof repeats the argument in [Renault 1980,  Appendix]  \cite{R}.
Corollary \ref{cor2} follows.

\begin{rmk}\label{rmk6}
\textnormal{
Corollary \ref{cor2} was first proved by  [Renault 1980,  Appendix]  \cite{R};
the GICAR algebra involved in the original proof  is isomorphic to the cluster $C^*$-algebra  ${\Bbb  A}(1,1)$, 
see Figure 3.
}
\end{rmk}

\subsection{The Hecke group}
The space  $T_{\goth A}$ can be parameterized by a Fuchsian group of the second
type.  Namely,  the {\it Hecke group} ${\goth H}(\lambda)$ is a subgroup of
$SL_2(\mathbf{R})$ generated by matrices:
\begin{equation}
\left(\matrix{0 & -1\cr 1 & 0}\right)\quad \hbox{and}
\quad \left(\matrix{1 & \lambda\cr 0 & 1}\right), 
\end{equation}
where $\lambda\in \mathbf{R}$. The group  ${\goth H}(\lambda)$
acts on the Lobachevsky plane ${\Bbb H}=\{z=x+iy ~|~ y>0\}$ by the linear fractional
transformations. The {\it Hecke theorem} says that  the action is discrete if and only if:    
$\lambda\in [2, \infty) \cup \{\lambda_n:=2\cos\left({\pi\over n}\right) ~|~ n\ge 3\}$. 
The fundamental domain $\{z\in {\Bbb H} ~|~ |\Re (z)|<{\lambda\over 2}, ~|z|>1\}$  of  ${\goth H}(\lambda)$ 
is sketched in Figure 4. It is easy to see, that ${\Bbb H}/{\goth H}(\lambda)$
is a disk if $\lambda\ge 2$ and there exits an Ahlfors map 
$f: {\goth A}\to {\Bbb H}/{\goth H}(\lambda)$ of degree $2$,  such that: 
\begin{equation}\label{eq35}
T_{\goth A}\cong \{t=\lambda^2 ~|~t\in [4, \infty) \cup \{4\cos^2 \left({\pi\over n}\right) ~|~ n\ge 3\}\}.  
\end{equation}
(In fact, the Hecke theorem follows from  (\ref{eq35}) and analysis of  algebra  ${\cal A}(1,1)$
given in Section 4.2;  the ${\cal A}(1,1)$ is a coordinate ring of  $T_{\goth A}$.)

\begin{figure}
\begin{picture}(300,100)(-40,0)


\qbezier(70,30)(90,70)(110,30)

\put(50,30){\line(1,0){80}}

\put(60,30){\line(0,1){40}}
\put(120,30){\line(0,1){40}}

\put(45,15){$-{\lambda\over 2}$}
\put(118,15){${\lambda\over 2}$}

\put(80,0){$\lambda\ge 2$}


\qbezier(220,30)(240,70)(260,30)

\put(200,30){\line(1,0){80}}

\put(230,45){\line(0,1){20}}
\put(250,45){\line(0,1){20}}

\put(210,15){$-1$}
\put(260,15){$1$}

\put(235,0){$\lambda_3$}

\end{picture}
\caption{The fundamental domain of the Hecke group.} 
\end{figure}
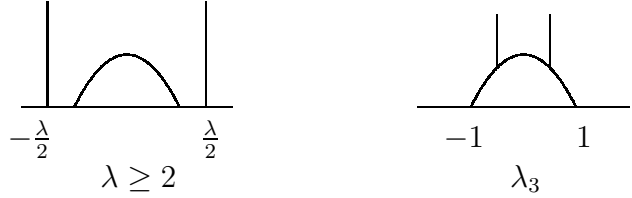

\begin{rmk}\label{rmk4}
\textnormal{
Formula (\ref{eq35}) relates the Hecke group  ${\goth H}(\lambda)$ to the algebra
${\Bbb A}(1,1)$ and henceforth to the index of subfactors by taking 
the weak closure of a representation of ${\Bbb A}(1,1)\subset M_{2^{\infty}}$;
such a question was raised in [Jones 1991, Section 3.1]  \cite{J1}.  
}
\end{rmk}

\bigskip\noindent
{\sf Acknowledgments.} 
I am grateful to  Ralf Schiffler for  an invitation to his Cluster Algebra Seminar at the University of Connecticut  and  many  helpful discussions.  
I thank the anonymous referee for the interest,  careful reading and useful suggestions regarding  this  paper.



\vskip1cm

\textsc{Department of Mathematics and Computer Science, St.~John's University, 8000 Utopia Parkway,  New York,  NY 11439, United States;}
~\textsc{E-mail:} {\sf igor.v.nikolaev@gmail.com}

\end{document}